\documentclass[12pt]{amsart}
\usepackage{txfonts}
\usepackage{amssymb}
\usepackage{eucal}
\usepackage{graphicx}
\usepackage{amsmath}
\usepackage{amscd}
\usepackage[all]{xy}           %xypic macro for latex2.09

\usepackage{amsfonts,latexsym}
\usepackage{xspace}
\usepackage{epsfig}
\usepackage{float}
\usepackage{color}
\usepackage{fancybox}
\usepackage{colordvi}
\usepackage{multicol}
\usepackage{colordvi}
\usepackage{stmaryrd}
\usepackage[colorlinks,final,backref=page,hyperindex,hypertex]{hyperref}

\newtheorem{prop-def}{Proposition-Definition}[section]

\newcommand{\nc}{\newcommand}
\newcommand{\delete}[1]{}

\nc{\mlabel}[1]{\label{#1}}  % Use this to suppress names
\nc{\mcite}[1]{\cite{#1}}  % Use this to suppress names
\nc{\mref}[1]{\ref{#1}}  % Use this to suppress names
\nc{\mbibitem}[1]{\bibitem{#1}} % Use this to show number
\delete{
\nc{\mlabel}[1]{\label{#1}  % Use the next two lines to show names
{\hfill \hspace{1cm}{\bf{{\ }\hfill(#1)}}}}
\nc{\mcite}[1]{\cite{#1}{{\bf{{\ }(#1)}}}}  % Use this lines to show names
\nc{\mref}[1]{\ref{#1}{{\bf{{\ }(#1)}}}}  % Use this lines to show names
\nc{\mbibitem}[1]{\bibitem[\bf #1]{#1}} % Use this to show name
}

\nc{\bfk}{\mathbf{k}}
\nc{\Der}{\mathrm{Der}}
\nc{\Ker}{\mathrm{Ker}}

\begin{document}

\title{Generalized derivations of $3$-Lie algebras }

\author{RuiPu  Bai}
\address{College of Mathematics and Information  Science,
Hebei University, Baoding 071002, China} \email{bairuipu@hbu.edu.cn}

\author{QiYong Li}
\address{College of Mathematics and Computer Science,
Hebei University, Baoding 071002, P.R. China}
\email{liqiyong2012@163.com}

\author{Kai Zhang}
\address{Library, Agricultural University of Hebei, Baoding, 071001, P. R. China}
\email{zhangkai@hebau.edu.cn}

\date{}

\subjclass[2000]{17B05, 17B30.}

\keywords{ 3-Lie algebra,  derivation, generalized derivation, quasiderivation, quasicentroid
}

\maketitle

\begin{abstract} Generalized derivations, quasiderivations and quasicentroid of $3$-algebras are introduced, and basic relations between them are studied. Structures of quasiderivations and quasicentroid of $3$-Lie algebras, which contains a maximal diagonalized tours, are systematically investigated.  The main results are:  for all $3$-Lie algebra $A$,  1)  the generalized derivation algebra $GDer(A)$ is the sum  of  quasiderivation algebra $QDer(A)$ and quasicentroid
$Q\Gamma(A)$; 2)  quasiderivations of $A$  can be embedded as derivations in a larger algebra;
3)  quasiderivation algebra $QDer(A)$ normalizer quasicentroid, that is, $[QDer(A), Q\Gamma(A)]\subseteq Q\Gamma(A)$; 4) if $A$ contains a  maximal diagonalized tours $T$, then $QDer(A)$ and $Q\Gamma(A)$ are diagonalized  $T$-modules, that is, as $T$-modules, $(T, T)$  semi-simplely acts on $QDer(A)$ and $Q\Gamma(A)$, respectively.

\end{abstract}

\baselineskip=18pt

\section{Introduction}
\setcounter{equation}{0}
\renewcommand{\theequation}
{1.\arabic{equation}}

  $3$-Lie algebras \cite{F} have close relationships with many
important fields in mathematics and mathematical physics (cf.
\cite{N, T, BL,HHM,HCK,G,P}). The multiple $M2$-brane model of Bagger-Lambert \cite{BL} and Gustavsson  \cite{G} is defined
on Lie 3-algebra, which serve as the gauge symmetry algebras for the $M2$-brane
world-volume theory. For the consistency of these symmetries, it need to impose the
fundamental identity on the Lie $3$-algberas.
P. Ho, Y. Imamura, Y. Matsuo in paper  \cite{HIM} studied two
derivations of the multiple $D2$ action from the multiple $M2$-brane
model proposed by Bagger-Lambert and Gustavsson. One defined $3$-Lie algebra which contains generators of a given Lie algebra. Such an extension contains generators with negative
norms. By suitably choosing such extension, one might restrict the
field associated with it to constant or zero while keeping almost all of the symmetry
of BLG theory.  The second derivation of multiple $D2$-brane, the extra generator has
a simple physical origin, the winding of $M5$-brane around $S^1$
which defines the reduction
from $M$-theory to the type IIA theory.
  One may provide a similar geometrical origin to
other $3$-Lie algebras.

The concept of derivations appear in different mathematical fields with many different forms.
In algebra systems,  derivations are linear mappings satisfying the Leibniz relation.
     So several kinds of  derivations, for example,  triple derivations, generalized derivations, qusiderivations  of algebras are studied \cite{WY, L, C, YZ}. In
      $n$-Lie algebras derivations are very useful, derivations are also studied \cite{B1, B2, B3, B4}. It is proved that the Lie algebra of the autoisomorphism group $Aut(L)$ of an $n$-Lie algebra $L$ is its derivation algebra $Der(L)$. And if $M$ is an $L$-module, then $M$ is also a  module of its inner derivation algebra ad$(L)$. By means of this relation,  finite dimensional irreducible modules  of simple $n$-Lie algebras over the complex field is classified \cite{D}.

In this paper, we study the  generalized derivations of $3$-Lie algebras.
Although  definitions of generalized derivations,  quasiderivations and  quasicentroid on $3$-Lie algebras are similar to the definitions
 on Lie algebras  \cite{L2},  but structures of them are very different ( see the discussion in section 4 and  section 5). In the following,  $\mathbb
F$ denotes  a field with   $ch \mathbb F\neq 2$ and $3$,  $Hom(W, V)$ is a vector space spanned by all the linear maps from vector space $W$  onto vector space $V$.

\section{Fundamental notions}
\setcounter{equation}{0}
\renewcommand{\theequation}
{2.\arabic{equation}}

{\it A $3$-Lie algebra} is a vector space $ A$ over $\mathbb
F$ endowed with a $3$-ary multi-linear skew-symmetric multiplication
$[ ~, ~, ~]$ satisfying following identity, for all $ x_1, x_2, x_3, y_2, y_3\in A,$
$$
  [[x_1, x_2, x_3], y_2, y_3]=\sum_{i=1}^3[x_1,  [ x_i, y_2, y_3],
  x_3].\eqno(2.1)$$

{\it A derivation } of a $3$-Lie algebra $A$ is
 an element $f\in $Hom$ (A, A)$ such that  for all  $x, y, z\in A$,
    $$ [f(x),y,z]+[x, f(y), z]+[x, y, f(z)]=f([x, y, z]).
$$
 The set  of all derivations of $A$ is  denoted by $Der(A)$, which   is a subalgebra of the general linear  Lie
algebra $gl(A)$.

A left multiplication of $A$ determined by $x_1, x_{2}\in
A$ is a linear map $\mbox{ad}(x_1, ~x_{2}):$ $
~A\rightarrow ~A$, defined by
\[
    \mbox{ad}(x_1, ~x_{2})(x)=[x_1, x_2, ~ x],
\]
where $x\in A$. A left multiplication is a derivation. The linear
combinations of left multiplications are called inner derivations;
the set of all inner derivations is denoted by $\mbox{ad}(A),$ which
is an ideal of $Der(A)$.

 Let $B_1, B_2, B_3$ be subspaces of $A$.
The notation $[B_1, B_3, B_3]$ denotes the subspace of $A$ spanned by vectors $[x_1, x_2, x_3]$, where $x_i\in B_i, 1\leq i\leq 3$.  A subspace $B$  is called  {\it a subalgebra $($an ideal$)$} if $B$ satisfies
$[B,B,B]\subseteq B$ ($[B, A, A]\subseteq B$). In particular, the subalgebra generated by
 vectors $[x_1, x_2, x_3]$ for all  $x_1, x_2, x_3\in A$ is
called the {\it derived algebra} of $A$, which is denoted by $A^{1}$.
If $A^{1}=0$, then $A$ is called an  {\it abelian algebra}.

{\it The center} of $A$ is $Z(A)=\{x\in A ~ | ~[x, A, A]=0\}$. It is clear that  $Z(A)$ is an abelian ideal of $A$.

A {\it centre derivation}   of $A$ (see \cite{B5}) is an element $f\in Hom(A,A)$ satisfying $f(A)\subseteq Z(A)$ and $f(A^{1})=0.$
  The set of all  centre derivations is denoted by $ZDer(A)$, which is an ideal of $Der(A)$.

{\bf Definition 2.1} {\it Let $A$ be a $3$-Lie algebra,  $f_1$ be a linear map of $A$. If there exist linear maps $f_2, f_3, f'$ of $A$ satisfy for all $ x, y, z $ $\in A, $
$$
     [f_1(x), y, z]+[x, f_2(y), z]+[x, y, f_3(z)]=f'([x, y, z]),  \eqno(2.2)
$$
then $f_1$ is called a generalized derivation of $A$. The set of all quaternions $(f_1, f_2, f_3, f')$ is denoted by $\Delta(A)$, and the set of all generalized derivations is denoted by $ GDer(A)$. In the case of $f_1=f_2=f_3=f$, the linear map $f$ is called a quasiderivation of $A$, and the set of all  quasiderivations is denoted by $QDer(A)$.}

By the above definition, for any $3$-Lie algebra $A$, we have
$$
     \mbox{ad}(A)\subseteq  Der(A)\subseteq QDer(A)\subseteq GDer(A).
$$
And by the skew-symmetry of the multiplication, if $(f_1, f_2, f_3, f')\in\Delta(A)$,  then
$$
     [f_2(x), y, z]+[x, f_1(y), z]+[x, y, f_3(z)]=f'([x, y, z]).\eqno(2.3)
$$
It follows that   $f_1, f_2, f_3\in GDer(A)$ and $( f_{i_1}, f_{i_2}, f_{i_3}, f')\in\Delta(A)$, where $(i_1, i_2, i_3)$ is arbitrary $3$-ary permutation.

We know  that the {\it  centroid}  {\cite{ B5} of a $3$-Lie algebra $A$ is an associative algebra $\Gamma(A)$ composed  of
linear maps $f$ satisfying for all  $x, y,z \in A$,
   \[ [f(x), y, z]=[x, f(y), z]=[x, y, f(z)]=f([x, y, z]).\eqno(2.4)
\]
We give the generalization of  $\Gamma(A)$ as follows.

{\bf Definition 2.2} {\it Let $A$ be a $3$-Lie algebra. The quasicentroid of $A$, is denoted by
$ Q\Gamma(A)$, is a vector space spanned by all elements  $f\in Hom(A, A)$ which satisfies  for all  $x, y,z \in A$,
\[
    [f(x), y, z]=[x, f(y), z]=[x, y, f(z)]. \eqno(2.5)
\]
}

\section{  Generalized derivations}

 In this section we study  the relations between $GDer(A)$, $QDer(A)$, $Der(A)$, $ Q\Gamma(A)$ and $\Gamma(A)$ of a $3$-Lie algebra $A$. And
 we  see that  quasiderivations of $A$ can be embedded as derivations in a larger algebra.

 {\bf Theorem 3.1}. {\it Let $A$ be a 3-Lie algebra. Then  $GDer(A)$, $QDer(A)$, $Q\Gamma(A)$ and $\Gamma(A)$ are subalgebras of the general linear Lie algebra $gl(A)$.
In addition, if $Z(A)=0,$ then $ Q\Gamma(A)$ and $\Gamma(A)$ are abelian $3$-Lie algebras.}

{\bf Proof} For arbitrary $f, g\in GDer(A)$ and  $(f, f', f'', f'''),$~ $(g, g', g'',g''')\in \Delta(A),$  by Eq.(2.2) and a direct computation we have
\[
[[f, g](x),y,z]+[x, [f', g'](y),z]+[x, y, [f'', g''](z)]=[f''', g'''][x, y, z], \text{for all} ~x, y, z\in A,
\]
then  $[f, g]=gf-fg\in GDer(A)$, and
$([f, g], [f', g'], [f'', g''], [f''', g'''])\in \Delta(A).$  It follows that $GDer(A)$ and $QDer(A)$ are subalgebras of $gl(A)$.

 For all $f, g\in Q\Gamma(A)$ and $x, y, z\in A$, we have

$[[f, g](x), y, z ]=[(fg-gf)(x), y, z)]=[g(x), f(y), z]-[f(x), g(y), z]$

$=[x, f(y), g(z)]-[f(x), y, g(z)]=[x, y, fg(z)]-[x, y, fg(z)]=0.$
\\
Since $Z(A)=0$, we obtain $[f, g]=0$.  $\Box$

 {\bf Lemma 3.2}. {\it Let $A$ be a 3-Lie algebra. Then  the following inclusions hold
$$Q\Gamma(A)\subseteq GDer(A),\eqno(3.1)$$
$$\Gamma(A)\subseteq QDer(A)\cap Q\Gamma(A),\eqno(3.2)$$
$$[Der(A), \Gamma(A)]\subseteq\Gamma(A),\eqno(3.3)$$
$$[QDer(A), Q\Gamma (A)]\subseteq Q\Gamma(A).\eqno(3.4)$$
}

{\bf Proof} For all $f\in Q\Gamma(A)$, we have
$$
[f(x), y, z]+[x, \frac{-1}{2}f(y), z]+[x, y, \frac{-1}{2}f(z)]=0,
$$
it implies that $(f, \frac{-1}{2}f, \frac{-1}{2}f, 0)\in \Delta(A)$. Therefore,  $f\in GDer(A)$.
Similarly, by direct computations, we obtain (3.2), (3.3) and (3.4). $\Box$

   The following examples to show  the relations between  $Der(A)$, $QDer(A)$ and $GDer(A)$.

  If $A$ is an abelian $3$-Lie algebra, then it is clear that  $Der(A)$ $=QDer(A)=GDer(A)=gl(A)$.

 If $A$ is a $3$-dimensional $3$-Lie algebra with $A^1\neq 0.$
Then there exists a basis $\{x_{1}, x_{2}, x_{3}\}$ such that  the multiplication is $[x_{1},$ $x_{2},x_{3}]=x_{1}$.
For all linear map $f$ of $A$, suppose $f(x_{i})=\sum\limits_{j=1}^{3}a_{ij}x_{j}$, $a_{ij}\in \mathbb F$. Define a linear map $f'$ of $A$ by
   $f'(x_{i})=\sum\limits_{j=1}^{3}b_{ij}x_{j},$ where  $b_{ij}\in \mathbb F$ satisfying
   $b_{11}=a_{11}+a_{22}+a_{33},~~ b_{12}=b_{13}=0.$ Then $(f, f, f, f')\in \Delta(A),$ that is, $f\in QDer(A)$. By a direct computation, $Der(A)\neq gl(A).$ Therefore, $Der(A)\varsubsetneq QDer(A)=gl(A).$

  Let $B$ be a $4$-dimensional $3$-Lie algebra with  the multiplication $[x_{1},x_{2},x_{3}]=x_{1}$, where  $\{x_{1}, x_{2}, x_{3}, x_4\}$ is a basis of $B$, then $Z(B)=\mathbb F x_4\neq 0$.
 For all  $(f, f, f, f') \in \Delta(A)$, suppose $f(x_i)=\sum\limits_{k=1}^4a_{ik}x_k$. Then $[x_i, x_j, f(x_4)]=[x_i, x_j, \sum\limits_{k=1}^4a_{4k}x_k]=0$.
We obtain $a_{41}=a_{42}=a_{43}=0.$  It follows $QDer(A)\neq gl(A)$. By a direct computation,   $Der(A)\neq QDer(A)$. Therefore,  $Der(A)\varsubsetneq QDer(A)\varsubsetneq gl(A).$

 {\bf Theorem 3.3}. {\it Let $A$ be a $3$-Lie algebra over  $\mathbb F$.  Then

(1) $GDer(A)=QDer(A)+Q\Gamma(A)$ (the sum of subspaces).

(2) $Q\Gamma(A)$ is an ideal of $GDer(A).$ In addition, if $Z(A)=0$, then $Q\Gamma(A)$ is an abelian ideal of $GDer(A)$.}

{\bf Proof} For every $g\in GDer(A)$, and   $(g,g',g'',g''')\in \Delta(A)$, we see that

\vspace{2mm}$(g,g',g'',g''')=(\frac{g+g'+g''}{3},\frac{g+g'+g''}{3},\frac{g+g'+g''}{3},g''')+(\frac{2g-g'-g''}{3},\frac{2g'-g-g''}{3},\frac{2g''-g-g'}{3},0)$.

Now  we prove  that $\frac{g+g'+g''}{3}\in QDer(A)$, and maps $\frac{2g-g'-g''}{3},\frac{2g'-g-g''}{3},\frac{2g''-g-g'}{3}$  $\in Q\Gamma(A)$.

From Eq.(2.2), we find that
\begin{align*}
&[\frac{g+g'+g''}{3}(x),y,z]+[x,\frac{g+g'+g''}{3}(y),z]+[x,y,\frac{g+g'+g''}{3}(z)]\\=&\frac{1}{3}([g(x),y,z]+[x,g(y),z]+[x,y,g(z)])
+\frac{1}{3}([g'(x),y,z]+[x,g'(y),z]+[x,y,g'(z)])\\+&\frac{1}{3}([g''(x),y,z)+[x,g''(y),z]+[x,y,g''(z)])\\=&g'''([z,y,z]).
\end{align*}
It follows that $\frac{g+g'+g''}{3}\in QDer(A)$.

By Eqs.(2.2) and (2.3), we have
\begin{align*}&[\frac{2g-g'-g''}{3}(x),y,z]-[x,\frac{2g-g'-g''}{3}(y),z]\\
=&\frac{2}{3}[g(x),y,z]-[\frac{1}{3}g'(x),y,z]-[\frac{1}{3}g''(x),y,z]
-[x,\frac{2}{3}g(y),z]+[x,\frac{1}{3}g'(y),z]+[x,\frac{1}{3}g''(y),z]\\
=&\frac{2}{3}g'''([x,y,z])-\frac{2}{3}[x,g'(y),z]-\frac{2}{3}[x,y,g''(z)]
-\frac{1}{3}g'''([x,y,z])+\frac{1}{3}[x,g(y),z]+\frac{1}{3}[x,y,g''(z)]\\
-&\frac{1}{3}g'''([x,y,z])+\frac{1}{3}[x,g'(y),z]+\frac{1}{3}[x,y,g(z)]
-\frac{2}{3}[x,g(y),z]+\frac{1}{3}[x,g'(y),z]+\frac{1}{3}[x,g''(y),z]\\
=&-\frac{1}{3}[x,y,g''(z)]+\frac{1}{3}[x,y,g(z)]-\frac{1}{3}[x,g(y),z]
+\frac{1}{3}[x,g''(y),z]\\
=&\frac{1}{3}([x,y,(g-g'')(z)]+[x,(g''-g)(y),z]),
\end{align*}

\hspace{3mm}$ g'''([x,y,z])=[g'(x),y,z]+[x,g(y),z]+[x,y,g''(z)], $

\hspace{3mm}$g'''([x,y,z])=[g'(x),y,z]+[x,g''(y),z]+[x,y,g(z)].$

Therefore,

\hspace{3mm}$[\frac{2g-g'-g''}{3}(x),y,z]=[x,\frac{2g-g'-g''}{3}(y),z]=[x,y,\frac{2g-g'-g''}{3}(z)],$

\hspace{3mm}$[\frac{2g'-g-g''}{3}(x),y,z]=[x,\frac{2g'-g-g''}{3}(y),z]=[x,y,\frac{2g'-g-g''}{3}(z)],$

\hspace{3mm}$[\frac{2g''-g-g'}{3}(x),y,z]=[x,\frac{2g''-g-g'}{3}(y),z]=[x,y,\frac{2g''-g-g'}{3}(z)].$

The result (1) follows.
The result (2) follows from Lemma 3.2. $\Box$

 If $A$ is a decomposable $3$-Lie algebra, that is, $A=H \oplus K$ (a direct sum of ideals). Then we have the following result.

 {\bf Theorem 3.4}. {\it Let  $A=H\oplus K$ and $Z(A)=0$. Then   for all  $g\in GDer(A)$, we have
 $g(H)\subseteq H$, $g(K)\subseteq K$, and
 $$GDer(A)=GDer(H)\oplus GDer(K), ~~QDer(A)=QDer(H)\oplus QDer(K).$$
}

{\bf Proof } The result follows from a direct computation. $\Box$

   \vspace{2mm} At last of this section, we show that quasiderivations of a $3$-Lie algebra  can be embedded as derivations in a larger algebra.
    For this, we need to construct $3$-Lie algebras by tensor product [20].

    Let $A$ be a $3$-Lie algebra over $\mathbb F$ and $t$ be an indeterminant.   Then quotient algebra $t\mathbb F[t]/(t^{4})$ is a $3$-dimensional commutative associative algebra.
     Let $\widetilde{A}$ be the  $3$-Lie algebra over $\mathbb F$ by the tensor product of $A$ and $t\mathbb F[t]/(t^{4})$, that is, $\widetilde{A}=A\otimes (t\mathbb F[t]/(t^{4}))$, and  write $at$, $at^{2}$ and $ at^{3}$ in place of $a\otimes t$, $a\otimes t^{2} $  and $a\otimes t^{3}$, respectively. Then the multiplication of $\widetilde{A}$ is
  $$[a_1t^l, a_2t^1, a_3t^1]=[a_1, a_2, a_3]t^{3}, \forall a_1, a_2, a_3\in A,\eqno(3.5)$$
  and $ [a_1t^l, a_2t^m, a_3t^n]=0 $ in the case that at least one of $l, m, n$ is larger than $1$.  Let $U$ be a subspace of $A$ satisfying $A=U\oplus A^{1}$, that is, the subspace $U$ is the complement of the derived algebra $A^{1}=[A, A, A]$, then
  $$\widetilde{A}=At\oplus At^{2}\oplus At^{3}=At\oplus At^{2}\oplus A^{1}t^{3}\oplus Ut^{3}. \eqno(3.6)$$

Define linear map $l_{u}:QDer(A)\longmapsto Hom(\widetilde{A}, \widetilde{A})$, by for all $(f,f,f,f')\in \Delta(A),$ and $a, b\in A, c\in A^{1}, u\in U$,
  $$l_{u}(f)(at+bt^{2}+ct^{3}+ut^{3})=f(a)t+f'(c)t^{3}.\eqno(3.7)$$
  \\  Then  $l_u$ is injective, and $l_{u}(f)$ does not depend on the choice of $f'$.

 {\bf Theorem 3.5.} {\it Let  $A$ be a $3$-Lie algebra,  $\widetilde{A}$ and $l_{u}$  be  defined as Eqs.$(3.6)$ and $(3.7)$, respectively.
Then $l_u(QDer(A))\subseteq Der(\widetilde{A})$. In addition, if  $Z(A)=0$, then $Der( \widetilde A)$ has a semidirect summation
$Der(\widetilde A)=l_{u}(QDer(A))\oplus Z(Der(\widetilde{A}))$.}

{\bf Proof }
  For all $a_{i}t+b_{i}t^{2}+c_{i}t^{3}+u_{i}t^{3}\in \widetilde{A}, i=1, 2, 3$,  by Eqs.(3.6) and (3.7),
    \begin{align*}
  &[l_{u}(f)(a_{1}t+b_{1}t^{2}+c_{1}t^{3}+u_{1}t^{3}), a_{2}t+b_{2}t^{2}+c_{2}t^{3}+u_{2}t^{3}, a_{3}t+b_{3}t^{2}+c_{3}t^{3}+u_{3}t^{3}]
  \\+&[a_{1}t+b_{1}t^{2}+c_{1}t^{3}+u_{1}t^{3},l_{u}(a_{2}t+b_{2}t^{2}+c_{2}t^{3}+u_{2}t^{3}),a_{3}t+b_{3}t^{2}+c_{3}t^{3}+u_{3}t^{3}]
  \\+&[a_{1}t+b_{1}t^{2}+c_{1}t^{3}+u_{1}t^{3},a_{2}t+b_{2}t^{2}+c_{2}t^{3}+u_{2}t^{3},l_{u}(f)(a_{3}t+b_{3}t^{2}+c_{3}t^{3}+u_{3}t^{3})]
  \\=&[f(a_{1}),a_{2},a_{3}]t^{3}+[a_{1},f(a_{2}),a_{3}]t^{3}+[a_{1},a_{2},f(a_{3})]t^{3}
  \\=&f'([a_{1},a_{2},a_{3}])t^{3},
  \end{align*}
  \begin{align*}
  &l_{u}(f)([a_{1}t+b_{1}t^{2}+c_{1}t^{3}+u_{1}t^{3},a_{2}t+b_{2}t^{2}+c_{2}t^{3}+u_{2}t^{3},a_{3}t+b_{3}t^{2}+c_{3}t^{3}+u_{3}t^{3}])
  \\=&l_{u}(f)([a_{1},a_{2},a_{3}]t^{3})
  =f'([a_{1},a_{2},a_{3}])t^{3}.
  \end{align*}
 Therefore,  $l_{u}(f)\in Der(\widetilde{A})$. It follows that $l_{u}(QDer(A))\subseteq Der(\widetilde{A})$.

Now if  $Z(A)=0$, then by Eq.(3.5)  the center of $3$-Lie algebra $\widetilde{A}$ is
$$
Z(\widetilde{A})=At^{2}+At^{3},\eqno(3.8)
$$ and  for all $D\in Der(\widetilde{A})$,  $D(Z(\widetilde{A}))\subseteq Z(\widetilde{A})$. Then  for all  linear map
$$f: At\oplus At^{2}\oplus Ut^{3}\longmapsto At^{2}\oplus At^{3},$$
$f$ can be extended  to an element of
$Z(Der(\widetilde{A}))$ by taking $f(A^{1}t^{3})=0$.

For every $g\in Der(\widetilde{A})$ and $a\in A$, suppose
$$g(at)=a't+b't^{2}+c't^{3}+u't^{3}, \mbox {where } a', b'\in A, c'\in A^{1}, u'\in U.\eqno(3.9)$$
Define linear map $f: \widetilde{A} \rightarrow \widetilde{A}$ by  for all $a, b\in A, c\in A^{1}, u\in U,$
$$ f(at)=b't^{2}+c't^{3}+u't^{3}, ~~ f(ct^{3})=0, ~~f(bt^{2}+ut^{3})=g(bt^{2}+ut^{3}).  $$
 By Eq.(3.8), $f\in Z(Der(\widetilde{A}))$, and
$$(g-f)(at)=a't,~ (g-f)(ct^{3})=c't^{3}. $$
From  Eq.(3.9), there exist  linear maps  $h, h'$ of $A$ satisfying  for all $a\in A, c\in A^1, u\in U$,
$$h(a)=a', ~ h'(c)=c', ~h'(u)=0,$$
then   $(h, h, h, h')\in \Delta(A)$, that is, $h\in QDer(A)$. Thanks to Eq.(3.7), $$g-f=l_{u}(h).$$
The proof is completed.  $\Box$

 $$  \textbf{4. QUASIDERIVATIONS\quad OF\quad 3-LIE\quad ALGEBRAS} $$

 In this section,   we first describe quasiderivations of $3$-Lie algebras by cohomology theory, and then we study quasiderivation algebras of a class of $3$-Lie algebras which contains a maximal diagonalized tours $T$. For convenience, for a $3$-Lie algebra $A$ and for all $x, y, z\in A$, we denote  the multiplication $[ x, y, z ]$ by $\mu(x, y, z)$.

Let $(A, \mu)$ be a $3$-Lie algebra. Denote the kernel of the multiplication $\mu$ by
  $$Ker(\mu)=\{v\otimes w\otimes u~|~\mu(v, w, u )=0, ~v, w, u\in A\}.\eqno(4.1)$$

  For every linear map $f$ of $A$, define a linear map $f^{*}: A\otimes A\otimes A\rightarrow A\otimes A\otimes A$  by for all $v, w, u\in A$,
 $$f^{*}(v\otimes w\otimes u)
  =f(v)\otimes w\otimes u+v\otimes f(w)\otimes u+v\otimes w\otimes f(u). \eqno(4.2)$$

 \vspace{2mm} {\bf Theorem 4.1} {\it Let $(A,\mu)$ be a $3$-Lie  algebra and $f\in Hom(A, A)$. Then $f\in QDer(A,\mu)$ if and only if
 $$f^{*}(Ker(\mu))\subseteq Ker(\mu),\eqno(4.3)$$
 that is, for all $v, w, u\in A$, if $\mu(v, w, u)=0$, then $\mu(f(v), w, u)+\mu(v, f(w), u)+\mu(v, w, f(u))=0$.}

  {\bf Proof} For $f\in QDer(A)$ and  $(f, f, f, f')\in \Delta(A)$, by Eqs.(4.1) and (4.2),  for all $v\otimes w\otimes u\in Ker(\mu)$,
$$\mu(f^{*}(v\otimes w\otimes u ))=\mu(f(v), w, u)+\mu(v, f(w), u)+\mu(v, w, f(u))=f'(\mu(v,w, u))=0.$$
Then  $f^{*}(Ker(\mu))\subseteq Ker(\mu)$.

 Conversely, let  $U$  be a complement of the derived algebra $A^1$ in $A$, that is, $A=A^1\oplus U$.  For $f\in Hom(A, A)$,  if $f$ satisfies inclusion (4.3).
  Define linear map $ f': A\rightarrow A$ by for all $u\in U$, $f'(u)=0,$ and for every $z=\sum\limits_{i=1}^m\mu(v_i, w_i, u_i)\in A^1$,
  $$f'(z)=\sum\limits_{i=1}^m\mu (f^*(v_i\otimes w_i\otimes u_i)).$$  Then $f'$ is well defined.

 In fact,  if $z=\sum\limits_{i=1}^m\mu(v_i, w_i, u_i)$ and $z=\sum\limits_{j=1}^l\mu(v'_j, w'_j, u'_j),$  where $v_i, w_i, u_i, v'_i, w'_i, u'_i\in A$,  $1\leq i\leq m, $ $1\leq j\leq l,$ then
 $$\sum\limits_{i=1}^m(v_i\otimes w_i\otimes u_i)-\sum\limits_{j=1}^l(v'_j\otimes w'_j\otimes u'_j)\in Ker \mu.$$
Thanks to  inclusion (4.3),  $\mu (f^*(\sum\limits_{i=1}^m(v_i\otimes w_i\otimes u_i)))-\mu (f^*(\sum\limits_{j=1}^l(v'_j\otimes w'_j\otimes u'_j)))=0.$

 Therefore,  for all $v, w, u\in A,$
 $$
 \mu(f(v), w, u)+\mu(v, f(w), u)+\mu(v, w, f(u))=\mu(f^{*}(v\otimes w\otimes u ))=f'(\mu(v,w, u)). $$
 It follows $f\in QDer(A).$ $\Box$

   In papers \cite{ AI, AI1, K}, authors introduced modules and cohomology theory of $3$-Lie algebras.  For describing quasiderivations, we first recall some definitions. Let $A$ be a
$3$-Lie algebra, $V$ be a vector space. If there
exists a linear map
 $\alpha: A\wedge A\rightarrow End(V)$ satisfying

\vspace{1mm} $ \alpha([x, y, z], w)=\alpha(y, z)\alpha(x,
w)+\alpha(z,x)\alpha(y, w)+\alpha(x, y)\alpha(z, w), $

$
\hspace{2cm} 0 = \alpha(z, w)\alpha(x, y)-\alpha(x, y)\alpha(z, w)+\alpha([x, y, z],
w)+\alpha(z, [x, y, w])
$

\vspace{1mm}
\noindent for all $x, y, z, w\in A$, then $(V, \alpha)$ is called  a
representation of $A$, or $V$ is {\it an $A$-module}. If $\alpha(x, y)=0$ for all $x, y\in A$, then $V$ is  a trivial module.
If $V=A$ and $\alpha(x, y)=\mbox{ad}(x, y)$ for all $x, y\in A$, then $A$ is called the adjoin module of $A$.

   A $3$-Lie algebra  $A$ as an adjoint module, we define

   $C^0(A, A)=Hom(A, A)$,\quad  $C^1(A, A)=Hom(A\otimes A\otimes A, A)$, \quad  $C^2(A, A)=H(A^{\otimes 5}, A)$,
   \\ and
    $\delta_0: C^0(A, A)\rightarrow C^1(A, A), $ \quad $\delta_1: C^1(A, A)\rightarrow C^2(A, A): $ for all $ x_i\in A$, $1\leq i\leq 5,$

 \vspace{2mm}$ \delta_0(f)(x_1, x_2, x_3)=f([x_1, x_2, x_3])-[f(x_1), x_2, x_3]-[x_1, f(x_2), x_3]-[x_1, x_2, f(x_3)],$ \hfill(4.4)

\vspace{2mm}$\delta_1(f)(x_1, x_2, x_3, x_4, x_5)=\sum\limits_{i=1}^3f(x_1, \cdots, [x_i, x_4, x_5], \cdots, x_3])-f([x_1, x_2, x_3], x_4, x_5)$\\
\vspace{1mm}\hspace{4.3cm}$+\sum\limits_{i=1}^3([x_1, \cdots, f(x_i, x_4, x_5), \cdots, x_3])-[f(x_1, x_2, x_3), x_4, x_5].$\hfill(4.5)

Denote $Ker \delta_i=Z^i(A, A)\subseteq C^i(A, A)$, $B^{i+1}(A, A)=$ $Im(\delta_i)\subseteq $ $C^{i+1}(A, A)$, $i=0, 1$.
 Thanks to  Eqs.(4.4) and (4.5),   $\delta_1\delta_0=0.$

 Let $\dot A$ be the trivial module of $A$ on
the underlying vector space of $A$, and

$C^0(A, \dot A)=Hom(A, \dot A)$, $C^1(A, \dot A)=Hom(A\otimes A\otimes A, \dot A)$, $C^2(A, \dot A)=H(A^{\otimes 5}, \dot A)$.

Define
$\dot{\delta}_0: C^0(A, \dot A)\rightarrow C^1(A, \dot A), ~~ \dot{\delta}_1: C^1(A, \dot A)\rightarrow C^2(A, \dot A)$ by for all $ x_i\in A, 1\leq i\leq 5,$
$$ \dot{\delta}_0(f)(x, y, z)=f([x, y, z]),\eqno(4.6)$$
 $$\dot{\delta}_1(f)(x_1, x_2, x_3, x_4, x_5)=\sum\limits_{i=1}^3f(x_1, [x_i, x_4, x_5],  x_3)-f([x_1, x_2, x_3], x_4, x_5), \eqno(4.7)$$
then  $\dot\delta_1 \dot\delta_0=0,$ and have the following result.

        {\bf Theorem 4.2} {\it
               Let $A$ be a $3$-Lie algebra, $f, f'$ be linear maps of $A$. Then

                \vspace{1mm} 1) $f\in QDer(A)$ if and only if $\delta_0(f)\in B^{1}(A, \dot{A})$. More specifically, $(f,f,f,f')\in \Delta(A)$ if and only if $\delta_0 (f)=\dot{\delta_0}(f-f').$

               \vspace{1mm}2) If $\dot{\delta_0}(f)\in Z^{1}(A, A)$, then for all $x, y, z, u, v\in A$,
               $$[f([x,y,z]),u,v]=[f([x,u,v]),y,z]+[x,f([y,u,v]),z]+[x,y,f([z,u,v])].\eqno(4.8)$$

               In particular, if $(f,f,f,f')\in \Delta(A)$, then}
               \\                $[(f-f')([x,y,z]),u,v]=[(f-f')([x,u,v]),y,z]+[x,(f-f')([y,u,v]),z]+[x,y,(f-f')([z,u,v])].$

                {\bf Proof}  If  $f\in QDer(A)$ and $(f,f,f,f')\in \Delta(A)$, by Eq.(4.4), we find that
                $$\delta_0(f)(x, y, z)=(f-f')([x, y, z])=\dot\delta_0(f-f')(x, y, z).$$
                 Therefore, $\delta_0(f)\in B^{1}(A, \dot{A}),$ and
                $\delta_0 f=\dot{\delta_0}(f-f').$

                Conversely, if   $\delta_0(f)\in B^{1}(A, \dot{A}),$ then there exists $f'\in C^0(A, \dot A)$ such that
                 $\delta_0(f)=\dot \delta_0 f'$, that is,  for all $x, y, z\in A,$

              $\delta_0(f)(x, y, z)=f([x, y, z])-[f(x), y, z]-[x, f(y), z]-[x, y, f(z)]=\delta_0 f'([x, y, z])=f'([x, y, z]).$
                 \\
                   It follows that
                 $(f'-f)([x, y, z])=[f(x), y, z]+[x, f(y), z]+[x, y, f(z)].$
                The result 1) holds.

               If $f\in Hom(A, A)$ satisfies $\dot{\delta_0} f\in Z^{1}(A, A)$, then $\delta_1\dot{\delta_0}( f )=0.$ Thanks to Eq.(4.5), for all $x, y, z, u, v\in A,$

$\delta_1\dot{\delta_0} (f)(x, y, z, u, v)$

$=[\dot{\delta_0}(f)(x,u,v),y,z]+[x,\dot{\delta_0}(f)(y,u,v),z]+[x,y,\dot{\delta_0}(f)(z,u,v)]-[\dot{\delta_0}(f)(x,y,z),u,v] $

$=[f([x,u,v]),y,z]+[x,f([y,u,v]),z]+[x,y,f([z,u,v])]-[f[x, y, z], u, v]=0.$
\\
It follows Eq.(4.8).

If  $(f, f, f, f')\in \Delta(A)$, then $f\in QDer(A)$. Thanks to result 1), $\delta_0(f)=\dot\delta_0(f-f').$
From $\dot{\delta_0}(f)\in Z^{1}(A, A)$ and $\delta_0(f)\in Z^1(A, A)$, we obtain that  $f-f'$ satisfies Eq.(4.8).  $\Box$

 In the following of this section, we study structures of qusiderivation algebras of a class of $3$-Lie algebras which containing  a maximal diagonalized tours.

              Let $A$ be a $3$-Lie algebra over an algebraically closed  field $\mathbb F$ with $ch F=0$, which contains a maximal diagonalized tours $T$, that is, $T$ is an abelian subalgebra, and $A$ has a root-subspace decomposition  associative $T$ as follows

              $ A=\sum\limits_{\gamma\in \Omega}A_{\gamma}, ~
              A_{\gamma}=\{ x\in A~ | ~\mbox{for all} ~ t_1, t_2\in T, \mbox{ad}(t_1, t_2)(x)=\gamma(t_1, t_2)x, \},$  $A_0=T$, \hfill(4.9)
             \\ where $\Omega\subseteq (T\wedge T)^*-\{0\}$ - the dual space of $T\wedge T.$
              If $A_{\alpha}\neq 0,$  $\alpha$ is called a weight of $A$ associative $T$.

              For example, let $A$ be a $4$-dimensional $3$-Lie algebra over the complex field $\mathbb F$ (see Lemm 3.1 in \cite{B6}). Then except cases $(b^1)$ and $(c^2)$, the subalgebra spanned by $e_3, e_{4}$ is a maximal diagonalized tours of $A$ .

Now let $A$ be a $3$-Lie algebra containing a maximal diagonalized tours $T$. For all $t_1, t_2\in T$, we define linear mapping $d(t_1, t_2): Hom(A, A)\rightarrow Hom(A, A)$ by for all $f\in Hom(A, A)$ and $x\in A,$

       $(d(t_1, t_2)f)(x)=[t_1, t_2, f(x)]-f([t_1, t_2, x])=(\mbox{ad}(t_1, t_2)f-f\mbox{ad}(t_1, t_2))(x).$ \hfill(4.10)

        For convenient,   sometimes we replace $d(t_1, t_2)f$ by $(t_1, t_2)f$, for all $t_1, t_2\in T$, $f\in Hom(A, A)$.

{\it In the following of this section, let $A$ be a $3$-Lie algebra over an algebraically closed  field $\mathbb F$ with $ch \mathbb F=0$, which contains a maximal diagonalized tours $T$, and $A_0=T.$}

         {\bf Lemma 4.3} {\it  Let $A$ be a $3$-Lie algebra over  $\mathbb F$. Then  for all  $t_1, t_2, t_3, t_4\in T$ and  $f\in QDer(A)$,

        1) $(t_1, t_2)f\in QDer(A)$, that is, $(T, T)QDer(A)\subseteq QDer(A)$;

        2)  $(t_{1}, t_{2})(t_{3}, t_{4})f-(t_{3}, t_{4})(t_{1}, t_{2})f=0;$

        3)   $((t_{1},t_{2})^{n}\cdot f)(x)=\sum \limits_{k=0}^{n}(-1)^{k+1}C_{n}^{k}\mbox{ad}^{n-k}(t_{1}, t_{2})f\mbox{ad}^{k}(t_{1},
 t_{2})(x).$\hfill(4.11)}

 {\bf Proof }
For arbitrary $f\in QDer(A)$,
      $t_1, t_2\in T$ and $ x, y, z\in A$, suppose $ (f, f, f, f')\in \Delta(A)$, by Eq.(4.10), we have
               \begin{align*}
               &[(t_{1}, t_{2})f(x), y, z]+[x, y, (t_{1}, t_{2})f(y), z]+[x, y, (t_{1}, t_{2})f(z)]\\
               =&[\mbox{ad}(t_{1}, t_{2})f(x), y, z]+[f(x), \mbox{ad}(t_{1}, t_{2})(y),z]+[f(x), y, \mbox{ad}(t_{1}, t_{2})(z)]\\
               +&[\mbox{ad}(t_{1}, t_{2})(x), f(y), z]+[x, \mbox{ad}(t_{1}, t_{2})f(y), z]+[x, f(y), \mbox{ad}(t_{1}, t_{2})(z)]\\
               +&[\mbox{ad}(t_{1}, t_{2})(x), y, f(z)]+[x, \mbox{ad}(t_{1}, t_{2})(y), f(z)]+[x, y, \mbox{ad}(t_{1}, t_{2})f(z)]\\
               -&([f(\mbox{ad}(t_{1}, t_{2})(x)), y, z]+[\mbox{ad}(t_{1}, t_{2})(x), f(y), z]+[\mbox{ad}(t_{1}, t_{2})(x), y, f(z)])\\
               -&([f(x), \mbox{ad}(t_{1}, t_{2})(y), z]+[x, f(\mbox{ad}(t_{1}, t_{2})(y)), z]+[x, \mbox{ad}(t_{1}, t_{2})(y), f(z)])\\
                -&([f(x), y, \mbox{ad}(t_{1}, t_{2})(z)]+[x, f(y), \mbox{ad}(t_{1}, t_{2})(z)]+[x, y, f(\mbox{ad}(t_{1}, t_{2})(z))])\\
              =&\mbox{ad}(t_{1}, t_{2})f'([x, y, z])-f'\mbox{ad}(t_{1}, t_{2})([x, y, z]).
               \end{align*}
It follows that $((t_1, t_2)f, (t_1, t_2)f, (t_1, t_2)f, \mbox{ad}(t_1, t_2)f'-f'\mbox{ad}(t_1, t_2))\in \Delta(A).$ The result 1) follows.

 The result 2) follows from Eq.(4.10) and a direct computation. Now we prove Eq.(4.11).

In the case $n=2$, we have
             \begin{align*}
              ((t_{1}, t_{2})^{2}f)(x)&=[t_{1}, t_{2}, ((t_{1}, t_{2})f)(x)]-((t_{1}, t_{2})f)([t_{1}, t_{2}, x])\\
                                         &=[t_{1}, t_{2}, [t_{1}, t_{2}, f(x)]-f([t_{1}, t_{2}, x])]\\
                                         &-[t_{1}, t_{2}, f([t_{1}, t_{2}, x])]+f([t_{1}, t_{2}, [t_{1}, t_{2},x]])\\
                                         &=ad^{2}(t_{1}, t_{2})f(x)-2ad(t_{1}, t_{2})fad(t_{1}, t_{2})(x)+fad^{2}(t_{1},
                                         t_{2})(x).
             \end{align*}

             Suppose identity (4.11) is true for the case $n$, then
             \begin{align*}
             ((t_{1}, t_{2})^{n+1}f)(x)&=((t_1, t_2)((t_1, t_2)^nf))(x)
             \\&=((t_1, t_2)(\sum \limits_{k=0}^{n}(-1)^{k+1}C_n^{k}\mbox{ad}^{n-k}(t_{1}, t_{2})f\mbox{ad}^{k}(t_{1}, t_{2})))(x)\\
                          &=(\sum \limits_{k=0}^{n}(-1)^{k+1}C_n^{k}\mbox{ad}^{n+1-k}(t_{1}, t_{2})f\mbox{ad}^{k}(t_{1}, t_{2}))(x)\\
                          &+(\sum \limits_{k=0}^{n}(-1)^{k+2}C_n^{k}\mbox{ad}^{n-k}(t_{1}, t_{2})f\mbox{ad}^{k+1}(t_{1}, t_{2}))(x)\\
                                            &=\sum \limits_{k=0}^{n}(-1)^{k+1}C_{n}^{k}\mbox{ad}^{n+1-k}(t_{1}, t_{2})f\mbox{ad}^{k}(t_{1}, t_{2})(x)\\
                                            &+\sum \limits_{l=1}^{n+1}(-1)^{l+1}C_{n}^{l+1}\mbox{ad}^{n+1-l}(t_{1}, t_{2})f\mbox{ad}^{l}(t_{1}, t_{2})(x)\\
                                            &=\sum \limits_{k=0}^{n+1}(-1)^{k+1}C_{n+1}^{k}\mbox{ad}^{n+1-k}(t_{1}, t_{2})f\mbox{ad}^{k}(t_{1},
                                            t_{2})(x).
             \end{align*}
The proof is completed. $\Box$

 From Lemma 4.3,  we can assume  the Fitting decomposition of $QDer(A)$  associative $T$ as
               $$QDer(A)=QDer(A)_{0}\bigoplus\sum\limits_{\alpha\in \Lambda}QDer(A)_{\alpha}, \eqno(4.12)$$
where $~ \Lambda\subseteq (T\wedge T)^*-\{0\},$ and

  $QDer(A)_{0}=\{f~| ~f\in QDer(A),  ~\mbox{for all}~t_1, t_2\in T , \mbox{there exists a positive intiger}~ m  $ such that $(t_1, t_2)^m
  f=0\},$

   $ QDer(A)_{\alpha}=\{f~| ~f\in QDer(A),$ for all $~t_1, t_2\in T,$ there exists a positive  integer $ m$  such that $ ((t_1, t_2)-\alpha(t_1, t_2)I_d)^m
  f=0\}.$

 For $\alpha\in \Lambda$, if $QDer(A)_{\alpha}\neq 0$, then $\alpha$ is called a
   weight of $T$.

{\bf Lemma 4.4 } {\it  Let $A$ be a $3$-Lie algebra.
      Then for all $\alpha\in \Lambda$  and $\gamma\in \Omega$,
     $$ QDer(A)_{\alpha}(A_{\gamma})\subseteq A_{\alpha+\gamma}.$$}

{\bf Proof } For all $f\in QDer(A)_{\alpha}, $ $x_{\gamma}\in
A_{\gamma}$ and $t_1, t_2\in T$, suppose $f(x_{\gamma})=\sum\limits_{\beta\in
(T\wedge T)^*}x_\beta$, where  $x_\beta \in A_\beta$. Then by Eqs.(4.9)-(4.11) and Lemma 4.3,  there exists an
integer $m$ such that
$$((t_1, t_2)-\alpha(t_1, t_2)I_d)^m
f=0.$$
 Since \begin{align*}
&(((t_1, t_2)-\alpha(t_1, t_2)I_d)^m
f)(x_{\gamma})\\
=&\sum\limits_{s=0}^m(-1)^{(m-s)}C_m^s\alpha(t_1,
t_2)^{m-s}((t_1, t_2)^sf)(x_{\gamma})\\
=&\sum\limits_{s=0}^m(-1)^{(m-s)}C_{m}^{s}\alpha(t_1,
t_2)^{m-s}\sum\limits_{k=0}^s(-1)^{k+1}C_{s}^k\gamma(t_1,
t_2)^{k}ad(t_1, t_2)^{s-k}f(x_{\gamma})\\
=&\sum\limits_{\beta\in (T\wedge
T)^*}\sum\limits_{s=0}^m(-1)^{(m-s)}C_{m}^{s}\alpha(t_1,
t_2)^{m-s}\sum\limits_{k=0}^s(-1)^{k+1}C_{s}^k\gamma(t_1,
t_2)^{k}\beta(t_1, t_2)^{s-k}x_{\beta}\\
=&\sum\limits_{\beta\in (T\wedge
T)^*}\sum\limits_{s=0}^m(-1)^{(m-s-1)}C_{m}^{s}\alpha(t_1,
t_2)^{m-s}(\beta-\gamma)(t_1, t_2)^sx_{\beta}\\
=&-\sum\limits_{\beta\in (T\wedge T)^*}(\beta-\gamma-\alpha)(t_1,
t_2)x_{\beta}\\
=&(\gamma+\alpha)(t_1,
t_2)f(x_\gamma)-(t_1, t_2)f(x_\gamma),\\
 \end{align*}
we have $f(x_{\gamma})\in
A_{\alpha+\gamma}$.~~  $\Box$

{\bf Theorem 4.5} {\it Let $A$ be a $3$-Lie algebra. Then the action of $T$ on $QDer(A)$
defined by Eq.(4.10) is diagonally, that is,
$$QDer(A)_{\alpha}=\{ f ~ |~ f\in QDer(A), ~\mbox{for all}~  t_1, t_2\in T, (t_1, t_2)f=\alpha(t_1,
t_2)f \}. \eqno(4.13)$$

Therefore, for all $f\in QDer(A)_0$ and $t_1, t_2\in T$, $(t_1, t_2)f=0$. And if $f\in QDer(A)_0\cap Der(A)$, then  for all $t_1, t_2\in T$ and $\gamma\in \Omega$,  $\gamma(f(t_1), t_2)+\gamma(t_1, f(t_2))=0.$}

{\bf Proof }  By Eq.(4.10) and Lemma 4.4, for  all $f\in QDer(A)_{\alpha}$,  $x_{\gamma}\in
A_\gamma$ and  $t_1, t_2\in T,$
 \begin{align*} &((t_1, t_2)f)(x_\gamma)\\
 =&[t_1, t_2, f(x_{\gamma})]-f([t_1, t_2,
x_{\gamma}])\\
=&(\alpha+\gamma)(t_1, t_2)f(x_{\gamma})-\gamma(t_1,
t_2)f(x_{\gamma})=\alpha(t_1, t_2)f(x_{\gamma}).
\end{align*} Thne identity (4.13) holds.

If $f\in Der(A)\cap QDer(A)_0$, then  for all $\gamma\in \Omega$, $x_{\gamma}\in
A_{\gamma}$ and $x_{\gamma}\neq 0$, by Lemma 4.4
 \begin{align*}&[f(t_1), t_2, x_{\gamma}]+[t_1, f(t_2), x_{\gamma}]+[t_1, t_2,
f(x_{\gamma})]\\
=&(\gamma(f(t_1), t_2)+\gamma(t_1, f(t_2))x_{\gamma}+\gamma(t_1,
t_2)f(x_{\gamma})\\
=&\gamma(t_1, t_2)f(x_{\gamma}).
\end{align*}\\

Therefore, $\gamma(f(t_1), t_2)+\gamma(t_1, f(t_2)=0$.
The result follows. $\Box$

At last of this section, we give an example to show that there does not exist inherent  relations between $QDer(A)_\alpha$ and the derivation algebra $Der(A)$.

Let $A$ be a $3$-Lie algebra with a basis $x_1, x_2, x_3$, and the multiplication is
$[x_1, x_2, x_3]=x_1.$ Then $T=\mathbb Fx_2+\mathbb F x_3$ is a maximal  diagonalized tours, and $A=T\dot+A_1$, where $A_1=\mathbb Fx_1$.
For all $f\in Hom(A, A)$, suppose $f(x_i)=\sum\limits_{j=1}^3a_{ij}x_j$, $1\leq i\leq 3.$ Then the matrix form of $f$ in the basis is
$\left( {\begin{array}{*{20}{c}}
{{a_{11}}}&a_{12}&a_{13}\\
{{a_{21}}}&{{a_{22}}}&{{a_{23}}}\\
{{a_{31}}}&{{a_{32}}}&{ {a_{33}}}
\end{array}} \right),$ $a_{ij}\in \mathbb F, 1\leq i, j\leq 3.$
By a direct computation, we have

\vspace{1mm}1).~ $QDer(A)=gl(A)$, and for all $f\in gl(A)$, $(f, f, f, f')\in \Delta(A)$ if and only if $$f'= \left( {\begin{array}{*{20}{c}}
{{a_{11}}+a_{22}+a_{33}}&b_{12}&b_{13}\\
{{b_{21}}}&{{b_{22}}}&{{b_{23}}}\\
{{b_{31}}}&{{b_{32}}}&b_{33}
\end{array}} \right), b_{ij}\in \mathbb F, 1\leq i, j\leq 3.$$

\vspace{1mm} 2). $QDer(A)=QDer{A}_0\oplus QDer{A}_1\oplus QDer{A}_{-1},$ where
\begin{align*} QDer(A)_0&=\Big{\{} \left( {\begin{array}{*{20}{c}}
{{a_{11}}}&0&0\\
0&{{a_{22}}}&{{a_{23}}}\\
0&{{a_{32}}}&a_{33}
\end{array}} \right)~ | ~ a_{ij}\in \mathbb F, 1\leq i, j\leq 3\Big{\}}.\\
QDer(A)_{-1}&=\Big{\{} \left( {\begin{array}{*{20}{c}}
0&0&0\\
a_{21}&0&0\\
a_{31}&0&0
\end{array}} \right)~ | ~ a_{j1}\in \mathbb F, j=2, 3\Big{\}}.
\\QDer(A)_{1}&=\Big{\{} \left( {\begin{array}{*{20}{c}}
0&a_{12}&a_{13}\\
0&0&0\\
0&0&0
\end{array}} \right)~ | ~ a_{1j}\in \mathbb F, j=2, 3\Big{\}}.
\end{align*}

3).  $Der(A)=\Big{\{} \left( {\begin{array}{*{20}{c}}
{{a_{11}}}&0&0\\
{{a_{21}}}&{{a_{22}}}&{{a_{23}}}\\
{{a_{31}}}&{{a_{32}}}&{ - {a_{22}}}
\end{array}} \right)~ | ~ a_{ij}\in \mathbb F, 1\leq i, j\leq 3\Big{\}}.$

\hspace{9mm}$ad(A)=\Big{\{} \left( {\begin{array}{*{20}{c}}
{{a_{11}}}&0&0\\
{{a_{21}}}&0&0\\
{{a_{31}}}&0&0
\end{array}} \right)~ | ~ a_{i}\in \mathbb F, 1\leq i\leq 3\Big{\}}.$

Then  $QDer(A)_{-1}\subseteq Der(A)$,  $QDer(A)_{1}\nsubseteq Der(A)$,  $QDer(A)_0\nsubseteq Der(A)$, $Der(A)\nsubseteq QDer(A)_0$.

$$  \textbf{5. QUASICENTROLDS\quad OF\quad 3-LIE\quad ALGEBRAS} $$

  From Theorem 3.3, general derivation algebra of a $3$-Lie algebra is generated by   quasiderivations and  quasicentroid. In this section, we study the structure of quasicentroid  of $3$-Lie algebras.
   From section 3,  the centroid  $\Gamma(A)$ is contained in $QDer(A)$, and $\Gamma(A)$, $Q\Gamma(A)$ are associative algebras. In addition, for the centerless $3$-Lie algebra  $A$,  $\Gamma(A)$ and  $Q\Gamma(A)$ are commutative.

 {\bf Lemma  5.1}  {\it Let $A$ be a $3$-Lie algebra, then

  1) $[\Gamma(A),  Q\Gamma(A)]\subseteq Hom(A, Z(A)); $

  2)  if $f\in \Gamma(A) $, then $Ker(f)$ and $Im(f)$ are ideals of $A$;

  3) if $A$ is  indecomposable and for $ f\in \Gamma(A)$ satisfying that  $x^{2}$ does not divide the minimal polynomial of $f$,
              then $f$ is invertible;

  4) if $A$ is indecomposable and the centroid $\Gamma(A)$ consists of semisimple elements, then $\Gamma(A)$ is a field.}

{\bf Proof} The result follows from a direct computation. $\Box$

   From Lemma 3.2, quasicentroid $Q\Gamma(A)$ is a  $3$-Lie algebra $A$-module  under the action defined by for all $f\in Q\Gamma(A)$ and  $x, y, z\in A, $
     $$((x,y)
    f)(z)=( \mbox{ad}(x,y)f-f \mbox{ad}(x,y))(z).\eqno(5.1)$$
 Therefore $Q\Gamma(A)$ is a  $T$-module, and we have the following identities.

 {\bf Lemma 5.2}  Let $A$ be a $3$-Lie algebra, and $ch \mathbb F\neq 2$. Then  for all $f\in Q\Gamma(A)$ and $ x,y\in A,$
\begin{align*}
  [x, f(x), y]&=0, \\
  \mbox{ad}(x,y)ad(f(x),y)&=\mbox{ad}(f(x), y) ad(x, y),\\
 ((x,y)f)(z)&=((y,z)f)(x)=((z,x)
    f)(y), \\
  \mbox{ad}^m(f(x),y)&=\mbox{ad}^m(x,y) f^{m},~ \mbox{for all positive integer}~  m,\\
   \mbox{ad}^{m+1}(x,y) f& =\mbox{ad}(f(x),y) \mbox{ad}^m(x,y), ~ \mbox{for all positive integer}~  m.
   \end{align*}

    {\bf Proof } For all $f\in Q\Gamma(A)$ and $ x, y,z\in A,$ since $[x,f(x), y]=[f(x),x,y]$ and $ch \mathbb F\neq 2$, we have  $[x,f(x),y]=0$.

     Follows from $[\mbox{ad}(x, y), \mbox{ad}(f(x),y)]=\mbox{ad}([x,y,f(x)],y)+\mbox{ad}(f(x),[x,y,y])=0,$  that $\mbox{ad}(x,y)$ and $\mbox{ad}(f(x),y)$ are commutative.

The identity $((x,y)f)(z)=((y,z)f)(x)=((z,x)
    f)(y)$ follows from  the skew-symmetry of the multiplication of $3$-Lie algebras, directly.

     For the case $m=2$

        $\mbox{ad}^2(f(x),y)(z)=ad(f(x),y)\mbox{ad}(f(x),y)(z)=[f(x),y,[f(x),y,z]] =[f(x),y[x,y,f(z)]]$

        \hspace{2.6cm}$
        =[[f(x),y,x],y,f(z)]+[x,[f(x),y,y],f(z)]$
$+[x,y,[f(x),y,f(z)]] $

\hspace{2.6cm}$=\mbox{ad}^2(x,y)f^{2}(z).     $

         Then for the case $m$

        $\mbox{ad}^{m}(f(x),y)(z)=\mbox{ad}(f(x),y)\mbox{ad}^{m-1}(f(x),y)(z) =\mbox{ad}(f(x),y)\mbox{ad}^{m-1}(x,y) f^{m-1}(z)$

        \hspace{3.1cm}$=\mbox{ad}^{m-1}(x,y)\mbox{ad}(f(x),y)f^{m-1}(z)=\mbox{ad}^{m-1}(x,y)[f(x),y,f^{m-1}(z)])$

       \hspace{3cm} $=\mbox{ad}^{m-1}(x,y)([x,y,f^{m}])(z)=\mbox{ad}^{m}(x,y)f^{m}(z).$

       Similar discussion, the last identity holds.    $\Box$

 Let $A$ be a $3$-Lie algebra over an algebraically closed  field $\mathbb F$ of characteristic zero, which  contains a maximal diagonalized tours $T$. By the discussion in section 4, $A$ has decomposition  (4.9) associative $T$. And $Q\Gamma(A)$ as a $T$-module in Eq.(5.1) has  Fitting decomposition
 $$Q\Gamma(A)=Q\Gamma(A)_0+ Q\Gamma(A)_1,~ \mbox{where} ~ Q\Gamma(A)_1=\sum\limits_{\alpha\in \Pi}Q\Gamma(A)_{\alpha},\eqno(5.2)$$
where $\Pi\subseteq (T\wedge T)^*-\{0\}$, and  $ Q\Gamma(A)_{\alpha}=\{ f\in Q\Gamma(A)~ |$ ~  for all $ t_1, t_2\in T, (t_1, t_2)f=((t_1, t_2)-\alpha(t_1, t_2)I_d)^mf=0,$ for some positive inetger $ m\}.$

 {\it In the following we suppose that $A$ is a $3$-Lie algebra over an algebraically closed  field $\mathbb F$ with $ch \mathbb F=0$,  which  contains a maximal diagonalized tours $T$ with $A_0=T$.  }

  By the completely similar discussion to section 4, we have the following result.

  {\bf Theorem 5.3} {\it Let $A$ be a $3$-Lie algebra. Then

1) as a  $T$-module, the decomposition (5.2) of $Q\Gamma(A)$ is diagonally, that is,
$$Q\Gamma(A)_{\alpha}=\{ f ~ |~ f\in Q\Gamma(A), ~~ (t_1, t_2)f=\alpha(t_1,
t_2)f, ~~ \mbox{for all} ~ t_1, t_2\in T \}; \eqno(5.3)$$

2) for all $ \alpha\in \Pi$ and $\gamma\in \Omega$, $Q\Gamma(A)_{\alpha}(A_{\gamma})\subseteq A_{\alpha+\gamma}. $ Therefore,
 $Q\Gamma(A)(T)\subseteq T$.}

{\bf Proof} The result 1) follows from the similar discussion  to Lemma 4.4 and Theorem 4.5.

 For all   $f\in Q\Gamma(A)$ and $x, y, z\in T $, by  Lemma 5.2
                     $$\mbox{ad}^{m+1}(x,y)f(z)=\mbox{ad}^m(x,y)ad(f(x),y)(z)=\mbox{ad}(f(x),y)\mbox{ad}^m(x,y)(z)=0.$$
Therefore, $f(T)\subseteq T$.  $\Box$

 {\bf Theorem 5.4} {\it Let $A$ be a $3$-Lie algebra. Then

  \vspace{1mm}  1)  $Q\Gamma(A)_{0}(A_{1})\subseteq A_{1}.$

\vspace{1mm}2)  $Q\Gamma(A)_{1}(T)=0.$

\vspace{1mm}3)  $Q\Gamma(A)_{1}(A_1)\subseteq Z_A(T)=\{x\in A ~ | ~ [x, T, A]=0\}$, therefore,  $Q\Gamma(A)_{1}(A)\subseteq Z_A(T).$

\vspace{1mm}4) For $\alpha, \beta\in \Omega$, if $\alpha+\beta\neq 0,$ then $(A_{\alpha}, A_{\beta})Q\Gamma(A)_{0}=0,$ and

\hspace{4mm}$(A_{\alpha}, A_{-\alpha}) Q\Gamma(A)_{0}\subseteq Q\Gamma(A)_0$, $(A_{\alpha}, A_{-\alpha}) Q\Gamma(A)_{0}(A_1)=0$.

\vspace{2mm} 5)  $((T, A_{1})Q\Gamma(A)_{0})(T)=0. $ }

               {\bf Proof } Define linear map $\sigma: Q\Gamma(A)\otimes A\rightarrow A$ by for all $f\in Q\Gamma(A)$ and $ z\in A$,
               $$\sigma(f\otimes z)=f(z).$$
Thanks to Eq.(5.1), $\sigma((x, y)f\otimes z+f\otimes \mbox{ad}(x, y)z)=\mbox{ad}(x, y)\sigma(f\otimes z),$ that is, $\sigma$ is a module homomorphism.
               Therefore,  $Q\Gamma(A)_{0}(A_1)\subseteq A_{1}$.  The result 1) holds.

                Thanks to Theorem 5.3, $Q\Gamma(A)_1(T)\subseteq T$. But from result 1),  $Q\Gamma(A)_{1}(T)\subseteq A_1$. Therefore,
$Q\Gamma(A)_{1}(T)\subseteq A_{1}\cap T=0$. The result 2) holds.

The result 3) follows from the result 2) and a simple computation.

              For all $\alpha, \beta\in \Omega$,  $x_{\alpha}\in A_{\alpha}, y_{\beta}\in A_{\beta}, z\in A$ and
              $ h_{1},h_{2}\in T,$ $ f_{0}\in Q\Gamma(A)_{0}$,

              If $\alpha+\beta\neq 0$, then

               $((h_{1},h_{2})(x_{\alpha},x_{\beta})f_{0})(z)$

               $=-((x_{\alpha},x_{\beta})(h_{1},h_{2})f_{0})(z)+([(h_{1},h_{2}),(x_{\alpha},x_{\beta})]f_{0})(z)$

               $=(\alpha+\beta)(h_{1},h_{2})((x_{\alpha},x_{\beta})f_{0})(z),$ we get
               $$(x_{\alpha}, x_{\beta})f_0\in Q\Gamma(A)_{\alpha+\beta}\subseteq Q\Gamma(L)_1.\eqno(5.4)$$
               For all $\delta\in \Omega$, by Theorem 5.3 and the result 3),
               $$(x_{\alpha}, x_{\beta})Q\Gamma(A)_0(A_{\delta})\subseteq Q\Gamma(A)_{\alpha+\beta}(A_{\delta})\subseteq A_{\alpha+\beta+\delta},\eqno(5.5)$$
               and $(x_{\alpha}, x_{\beta})Q\Gamma(A)_0(A_{\delta})\subseteq T.$
               We find that
            if $ \alpha+\beta+\delta\neq 0$, then  $ (x_{\alpha}, x_{\beta})Q\Gamma(A)_0(A_{\delta})=0. $

                Now suppose $((x_\alpha, x_{\beta})Q\Gamma(A)_{0})(A_1)\neq 0$.
then there exist $\delta, \lambda, \mu\in \Omega$, and nonzero vectors $z_{\delta}\in A_{\delta},  m_{\lambda}\in A_{\lambda}, n_{\mu}\in A_{\mu},$ and $ f\in Q\Gamma(A)_{0}$
                such that $$[((x_{\alpha}, y_{\beta})f)(z_{\delta}), m_{\lambda},n_{\mu}]\neq 0.$$ Then
                $$ [((x_\alpha, x_{\beta})f)(z_{\delta}),  m_{\lambda},n_{\mu}]=[z_{\delta}, ((x_\alpha, x_{\beta})f)(m_{\lambda}) ,n_{\mu}]=
                [z_{\delta},  m_{\lambda},((x_\alpha, x_{\beta})f)(n_{\mu})]\neq 0.$$
                From above discussion, we have
$$ \alpha+\beta+\delta=0, ~~ \alpha+\beta+\lambda=0, ~~\alpha+\beta+\mu=0,
                \lambda+\mu=\delta+\mu=\delta+\lambda,$$
                and
                $\delta=\lambda=\mu=0$. Contradiction. It follows that if $\alpha+\beta\neq 0$, then $(A_\alpha, A_{\beta})Q\Gamma(A)_0=0.$

If $\alpha+\beta=0$, then for all $h_1, h_2\in T$ and $f_0\in Q\Gamma(A)_0$,
             $$(h_{1},h_{2})(x_{\alpha},x_{-\alpha}) f_{0}=(x_{\alpha},x_{-\alpha})(h_{1},h_{2})f_{0}+[(h_{1},h_{2}),(x_{\alpha},x_{-\alpha})]f_0=0.$$
              It implies that  $(A_{\alpha}, A_{-\alpha}) Q\Gamma(A)_0\subseteq  Q\Gamma(A)_0.$ Thanks to Lemma 5.2 and 5.3, for all $f_0\in \Gamma(A)_0,$ $\delta\in \Omega$,  $x_{\delta}\in A_{\delta}$,
               $$ (x_{\alpha}, x_{-\alpha}) f_0(x_{\delta})=(x_{-\alpha}, x_{\delta}) f_0(x_{\alpha})=(x_{\delta}, x_{\alpha}) f_0(x_{-\alpha})=0.$$
               Therefore,  $((x_{\alpha}, x_{-\alpha}) f_0)(A_1)=0$.
                              The result 4) follows.

               From Lemma  5.2 and Lemma 5.3,

              $(T,T)Q\Gamma(L)_{0}=0$, $ ((T, A_{1})Q\Gamma(A)_{0})(T)=((T,T)Q\Gamma(A)_{0})(A_1)=0.$  It follows result 5). The proof is completed. $\Box$

   {\bf Theorem 5.5} {\it
              Le $A$ be a $3$-Lie algebra with trivial centre, then $Q\Gamma(A)=\Gamma(A)\oplus Q\Gamma(A)_1$ with  $Q\Gamma(A)_1 Q\Gamma(A)_1=0$.}

             {\bf Proof } First, we prove that  $\Gamma(A)= Q\Gamma(A)_0$. It is clear that $\Gamma(A)\subseteq Q\Gamma(A)_0$. And from Lemma 5.2, for all $f_0\in Q\Gamma(A)_0$ and $ t_1, t_2\in T,$
             $(t_1, t_2)f_0=0$. Then for all $x\in A$ and  $t\in T$,
                 $$[t_1, t_2, f_0(x)]=f_0([t_1, t_2, x]),~~ [t_1, t_2, f_0(t)]=f_0([t_1, t_2, t])=0.$$

             Thanks to Theorem 5.4, for all $ \alpha, \beta\in \Omega$, $x_{\alpha}\in A_{\alpha}, ~x_{\beta}\in A_{\beta}, ~x\in A$ and $ f_0\in Q\Gamma(A)_0 $, if $\alpha+\beta\neq 0,$
            $((x_{\alpha}, $ $x_{\beta})f_0)(x)=0$. Then
             $$f_{0}([x_{\alpha}, x_{\beta}, x])=[x_{\alpha}, x_{\beta}, f_0(x)]-((x_{\alpha}, x_{\beta})f_0)(x)=[x_{\alpha}, x_{\beta}, f_0(x)].$$

 If $\alpha+\beta= 0,$ then by Theorem 5.4, for all $t\in T$, $\delta\in \Omega, x_{\delta}\in A_{\delta}, $  $((x_{\alpha}, x_{-\alpha})f_0)(x_{\delta})=0$.  Then
  $$f_{0}([x_{\alpha}, x_{-\alpha}, x_{\delta}])=[x_{\alpha}, x_{-\alpha}, f_0(x_{\delta})]-((x_{\alpha}, x_{-\alpha})f_0)(x_{\delta})=[x_{\alpha}, x_{-\alpha}, f_0(x_{\delta})].$$

 Since $((x_{\alpha}, x_{-\alpha})f_0)(t)\in T$,
 $$[((x_{\alpha}, x_{-\alpha})f_0)(t), T, T]=0, ~~ [((x_{\alpha}, x_{-\alpha})f_0)(t), x_{\delta}, A]=[t, ((x_{\alpha}, x_{-\alpha})f_0)(x_{\delta}), A]=0.$$
 It follows that $((x_{\alpha}, x_{-\alpha})f_0)(t)\in Z(A)=0$,  $((x_{\alpha}, x_{-\alpha})f_0)(T)=0.$ Therefore,
  $$f_{0}([x_{\alpha}, x_{-\alpha}, t])=[x_{\alpha}, x_{-\alpha}, f_0(t)]-((x_{\alpha}, x_{-\alpha})f_0)(t)=[x_{\alpha}, x_{-\alpha}, f_0(t)].$$

Summarizing above discussion, we have $f_0\in \Gamma(A)$, and   $\Gamma(A)=Q\Gamma(A)_0$.

  Thanks to Theorem 5.4,  $Q\Gamma(A)_1Q\Gamma(A)_1(A)\subseteq $ $Q\Gamma(A)_1(T)=0$. $\Box$

 {\bf Theorem 5.6} {\it Let $A$ be a $3$-Lie algebra with the decomposition $A=A_1\oplus A_2$, $[A_1, A_2, A]=0.$
Then we have
 $$Q\Gamma(A)=Q\Gamma(A_1)+Q\Gamma(A_2)+\Gamma_1+\Gamma_2,$$
 where $\Gamma_i=\{f\in Hom(A_i, A_j) ~|~f(A_i)\subseteq Z(A_j), 1\leq i\neq j\leq 2\}.$}

{\bf Proof }
It is clear that $Q\Gamma(A_1)+Q\Gamma(A_2)+\Gamma_1+\Gamma_2\subseteq Q\Gamma(A).$

Let $p_i\in Hom(A, A_i)$, $i=1, 2$ be projections, that is,  for all $x=x_1+x_2\in A$,  $ x_i\in A_i$, $p_i(x)=x_i$, $i=1, 2$.

For all $f\in Q\Gamma(A)$, define $f_i\in Hom (A_i, A_i)$, $i=1, 2$, $f_3\in Hom (A_1, A_2)$, $f_4\in Hom (A_2, A_1)$ as follows
$$f_1(x)=p_1f(x), ~\mbox{for all}~ x\in A_1; ~f_2(x)=p_2f(x), ~\mbox{for all}~ x\in A_2;$$
$$f_3(x)=p_2f(x), ~\mbox{for all}~  x\in A_1; ~f_4(x)=p_1f(x), ~\mbox{for all}~ x\in A_2,$$
 Without loss of generality, suppose $f_1(A_2)=f_3(A_2)=0$, $f_2(A_1)=f_4(A_1)=0$. Then for all $x=x_1+x_2\in A$, $x_1\in A_1$ and $x_2\in A_2$,
 $$f(x)=f(x_1)+f(x_2)=p_1f(x_1)+p_2f(x_1)+p_1f(x_2)+p_2f(x_2)=(f_1+f_2+f_3+f_4)(x).$$

Thanks to $[A, A_1, A_2]=0$, $f_1\in QC(A_1)$, $f_2\in QC(A_2)$.

For all $x\in A_1, y, z\in A_2$,
since

$[f_3(x), y, z]=$  $[f(x), y, z]=[x, f(y), z]=0,$
$[f_3(A_1), A_2, A_2]=0$. It follows that $f_3\in \Gamma_1$.

Similarly, we have $f_4\in \Gamma_2$. Therefore,  $Q\Gamma(A)\subseteq Q\Gamma(A_1)+Q\Gamma(A_2)+\Gamma_1+\Gamma_2.$  The result follows. $\Box$

 {\bf Corollary 5.7} {\it Let $A$ be a $3$-Lie algebra  with the decomposition $A=A_1\oplus\cdots\oplus A_m$, and  $[A_i, A_j, A]=0,$ $ 1\leq i\neq j\leq m.$
Then
$$Q\Gamma(A)=Q\Gamma(A_1)+\cdots +Q\Gamma(A_m)+\sum\limits_{1\leq i\neq j\leq m}\Gamma_{ij},$$
  where $\Gamma_{ij}=\{f\in Hom(A_i, A_j) ~|~f(A_i)\subseteq Z(A_j), 1\leq i\neq j\leq m\}.$}

{\bf Proof } The discussion is completely similar to Theorem 5.6. $\Box$

\section*{Acknowledgements}

The first author (R.-P. Bai) was supported in part by the Natural
Science Foundation (11371245) and the Natural
Science Foundation of Hebei Province, China (A2014201006).

\bibliography{}

\end{document}